\documentclass[letterpaper, 10pt, journal]{ieeetran}

\IEEEoverridecommandlockouts                              

\usepackage{cite}
\usepackage{graphicx}      

\graphicspath{Figures/}
\DeclareGraphicsExtensions{.png,PNG,.pdf,.PDF}
\usepackage{bbm}
\usepackage{comment}
\usepackage{bm}
\usepackage{amsmath} 
\protected\edef\ell{\noexpand\ensuremath{{\mathchar\the\ell}}}
\usepackage{amssymb}  
\usepackage{mathrsfs}
\usepackage{multicol}
\usepackage{xcolor}
\usepackage{mathtools}
\usepackage{soul}
\usepackage[T1]{fontenc}
\usepackage{eucal}
\usepackage{url}
\usepackage{amsmath}

\newtheorem{remark}{Remark}
\newtheorem{assumption}{Assumption}

\newtheorem{proposition}{Proposition}
\newtheorem{proof}{Proof}
\newtheorem{problem}{Problem}

\newtheorem{theorem}{Theorem}

\DeclareMathAlphabet{\pazocal}{OMS}{zplm}{m}{n}

\DeclareMathOperator*{\argmin}{arg\,min}

\title{ Homothetic tube model predictive control with multi-step predictors - extended version}

\author{Danilo Saccani, Giancarlo Ferrari-Trecate, Melanie N. Zeilinger, Johannes K{\"o}hler
\thanks{This research has been supported by the Swiss National Science Foundation under the NCCR Automation (grant agreement 51NF40\_180545) and Johannes K{\"o}hler was also supported by an ETH Career Seed Award funded
through the ETH Zurich Foundation.}
\thanks{D. Saccani and G. Ferrari-Trecate are with the Institute of Mechanical Engineering, Ecole Polytechnique Fédérale de Lausanne (EPFL), CH-1015 Lausanne, Switzerland. (email:  \{danilo.saccani, giancarlo.ferraritrecate\}@epfl.ch)  }%
\thanks{J. K{\"o}hler and M.N. Zeilinger are with the Institute for Dynamic Systems and Control, ETH Zurich, Switzerland. (email \{jkoehle,mzeilinger\}@ethz.ch)  }%
}

\begin{document}

\maketitle
\IEEEoverridecommandlockouts
\IEEEpubid{\begin{minipage}{\textwidth}\ \\[22pt] \\ \\
         \copyright 2023 IEEE.  Personal use of this material is  permitted.  Permission from IEEE must be obtained for all other uses, in  any current or future media, including reprinting/republishing this material for advertising or promotional purposes, creating new  collective works, for resale or redistribution to servers or lists, or  reuse of any copyrighted component of this work in other works.
     \end{minipage}}
\begin{abstract}                
We present a robust model predictive control (MPC) framework for linear systems facing bounded parametric uncertainty and bounded disturbances. 
Our approach deviates from standard MPC formulations by integrating multi-step predictors, which provide reduced error bounds. These bounds, derived from multi-step predictors, are utilized in a homothetic tube formulation to mitigate conservatism.
Lastly, a multi-rate formulation is adopted to handle the incompatibilities of multi-step predictors. 
We provide a theoretical analysis, guaranteeing robust recursive feasibility, constraint satisfaction, and (practical) stability of the desired setpoint.
We use a simulation example to compare it to existing literature and demonstrate advantages in terms of conservatism and computational complexity.
\end{abstract}



\section{Introduction}
The control of uncertain systems subject to safety critical constraints is a crucial problem in various engineering applications. To address this challenge, Model Predictive Control (MPC) has emerged as a powerful approach that can handle constraints efficiently~\cite{kouvaritakis2016model}. 
To ensure stability and robust constraint satisfaction, MPC requires a prediction model.
Such a model is obtained through system identification techniques, such as Bayesian linear regression~\cite{bishop2006pattern} or set-membership estimation~\cite{milanese1991optimal}, which also provide bounds on the parametric uncertainty. 
In this work, we propose a robust MPC formulation that can efficiently incorporate the uncertainty bounds from multi-step predictors.
\subsubsection*{Related work}
Robust MPC formulations ensuring recursive feasibility and stability for systems subject to bounded disturbances are well-established~\cite{kouvaritakis2016model}. 
In recent years, significant advances have been made to develop robust MPC methods that can effectively incorporate uncertainty bounds on the parameters of the state space model~\cite{langson2004robust,parsi2022scalable,lorenzen2019robust,Koehler2019Adaptive, Lu2019RAMPC,schwenkel2022model,chen2022robust}. 
These methods ensure safe operation by predicting a tube that bounds all possible uncertain systems. 
By using some nominal trajectory in combination with a scaled set, a simple \textit{homothetic} parametrization can be obtained~\cite{langson2004robust,lorenzen2019robust,Koehler2019Adaptive,parsi2022scalable}. 
By online adapting the parametric model, such robust approaches can be naturally extended to a robust \textit{adaptive} MPC formulations~\cite{lorenzen2019robust,Koehler2019Adaptive,Lu2019RAMPC}.

In recent years, an alternative paradigm based on multi-step predictors instead of traditional state space models is gaining traction~\cite{kohler2022state}. 
Multi-step predictors 
alleviate the need for complex disturbance propagation techniques. 
Early works in this direction focused on the special case of impulse response models, with a non-conservative treatment of parametric uncertainty~\cite{zheng1993robust,tanaskovic2014adaptive,furieri2022near}. 
Similar results have been presented with orthonormal basis functions~\cite{heirung2017dual}.
Recent developments provide identification routines with error bounds for general linear systems with measurement noise and disturbances~\cite{terzi2019learning, lauricella2020set}.
These methods directly yield prediction error bounds over multiple time steps, avoiding the need for conservative over-approximations~\cite{kohler2022state} and demonstrating their potential in experimental applications~\cite{lauricella2023day,saccani2022multitrajectory}.
For nonlinear and non-parametric models, similar multi-step predictors have been investigated~\cite{maddalena2021kpc, pfefferkorn2022exact, hashimoto2022learning,park2023simultaneous,he2022finite}.

Despite the availability of good error bounds for multi-step predictors, there is a gap in exploiting their full potential within robust MPC formulations. In~\cite{terzi2022robust,terzi2019learning} the prediction error due to the parametric uncertainty is over-approximated with a worst-case constant bound.
This leads to a simple MPC scheme, however, also introduces unnecessary conservatism. 
On the other hand, robust MPC approaches designed for state space models, e.g.,~\cite{langson2004robust,parsi2022scalable,lorenzen2019robust,Koehler2019Adaptive, Lu2019RAMPC,schwenkel2022model,chen2022robust}, can directly account for parametric uncertainty, thereby avoiding conservatism.
\subsubsection*{Contribution}
The main contribution of this paper is the development of a novel robust MPC formulation that fully utilizes the error bounds provided by multi-step predictors. The proposed formulation incorporates two key elements: (i) multi-rate scheme~\cite{terzi2022robust} for multi-step predictors and (ii) the homothetic tube formulation~\cite{lorenzen2019robust,Koehler2019Adaptive} to account for parametric uncertainty.
The homothetic tube formulation addresses parametric uncertainty in a non-conservative fashion, reducing the conservatism in~\cite{terzi2022robust}, where parametric uncertainty is conservatively upper-bounded within a rigid tube-based framework.
Finally, we provide a theoretical analysis that ensures robust constraint satisfaction and practical asymptotic stability.
We provide a qualitative discussion and a simulation example to contrast the proposed approach in terms of computational complexity and conservatism with state-of-the-art approaches~\cite{lorenzen2019robust,terzi2022robust,Koehler2019Adaptive}. This comparison demonstrates the advantages of the proposed approach: a reduced level of conservatism with moderate computational complexity.  
\subsubsection*{Notation}
We denote with $\mathbb{N}$ the set of non-negative integers and $\mathbb{N}_a^b=\{ n \in \mathbb{N} \ | \ a\leq n \leq b \}$. $[A]_i$ and $[a]_j$ denote the $j$-th row and entry of the
matrix $A$ and vector $a$, respectively. By 
$\bm{1}$ we denote a column vector of ones of appropriate size. 
Positive semidefinite matrices $A$ are denoted $A \succeq 0$. The quadratic norm of a vector $v$ w.r.t. a positive definite matrix $A$ is denoted by $\|v\|_A = \sqrt{v^\top A v}$. 
A function $\gamma:\mathbb{R}_{\geq0}\rightarrow\mathbb{R}_{\geq0}$ is of class $\mathcal{K}_\infty$ if it is continuous, strictly increasing, $\gamma(0) = 0$ and $\gamma(s)\rightarrow +\infty$ as $s\rightarrow +\infty$.
Given a set of points $\pazocal{S}^{v}=\{x^1,\ldots,x^{n_v}\}$, we use $\text{co}(\pazocal{S}^{v})$ to denote their convex hull. Lastly, for a sequence of matrices $M=\{ M_1, \dots, M_n \}$, we use $\text{blkdiag}(M)$ to denote the block diagonal matrix with $M_i$ placed along its diagonal.
\section{Problem formulation}
In this paper, we consider the following linear time-invariant discrete-time system with state $x_k\in\mathbb{R}^{n_x}$, input $u_k\in\mathbb{R}^{n_u}$ and additive bounded disturbance $w_k$:
\begin{equation} \label{eq:system}
    x_{k+1} = A x_k + B u_k + Mw_k,
\end{equation}
where $w_k\in\mathbb{W}\subseteq \mathbb{R}^{n_w}$ and the system matrices $(A,B)$ are not assumed to be known.
The state and input variables are constrained to lie within the following compact polytopes $\forall k \in \mathbb{N}$:
\begin{equation} \label{eq:constraints}
    x_k\in \mathbb{X}, \ \ \ \  u_k \in \mathbb{U},
\end{equation}
where $\mathbb{X}=\{ x\in \mathbb{R}^{n_x} | Fx\leq \bm{1} \}$ and $\mathbb{U}=\{ u\in \mathbb{R}^{n_u} | Gu\leq \bm{1} \}$ with given matrices $F\in\mathbb{R}^{c_{\mathrm{x}}\times n_x}$ and $G\in\mathbb{R}^{c_{\mathrm{u}}\times n_u}$.\\
For a given $p\in\mathbb{N}$, system~\eqref{eq:system} can be expressed in the following $p$-steps ahead state space representation:
\begin{align} \label{eq:multistepmodel}
    X_{j+1} &= \bar{A}(\theta)  X_j + \bar{B}(\theta)U_j  + \bar{M} W_j \nonumber \\
    Y_j &= \bar{C} (\theta) X_j+\bar{D} (\theta) U_j + \bar{N} W_j,
\end{align}
with state $X_j=x_{jp}$, input $U_j=[ u_{jp}^{\top},\ \dots, \ u_{jp+p-1}^{\top}]^{\top}\in\mathbb{U}^p$, additive disturbance $W_j=[w_{jp}^{\top}, \dots, w_{jp+p-1}^{\top}]^{\top}\in\mathbb{W}^p$, output $Y_j = [x_{jp+1}^{\top}, \dots, x_{jp+p-1}^{\top}]^{\top}\in\mathbb{X}^{p-1}$. 
Matrices $\bar{A}$, $\bar{B}$, $\bar{M}$, $\bar{C}$, $\bar{D}$ and $\bar{N}$ directly map the current system state to the future state sequence, instead of sequentially applying the one-step prediction of the state space model~\eqref{eq:system}.
We do not assume prior knowledge of the system matrices $(A,B)$ and hence the matrices in~\eqref{eq:multistepmodel} depend on some parameter vector $\theta$, which is unknown but constant. We make the following assumption about the disturbance and uncertainty.\begin{assumption}\label{ass:bound}
    The matrices in~\eqref{eq:multistepmodel} depend affinely on the parameter vector $\theta$
    \begin{multline}
        (\bar{A}(\theta),\bar{B}(\theta),\bar{C}(\theta),\bar{D}(\theta)) = (\bar{A}_0,\bar{B}_0,\bar{C}_0,\bar{D}_0)\\ +\sum_{i=1}^{n_{p}} (\bar{A}_i, \bar{B}_i, \bar{C}_i, \bar{D}_i) [\theta]_i.
    \end{multline}
    The true parameter vector $\theta^*$ lies inside a known bounded polytopic set:
    \begin{equation} \label{eq:modelUncer}
        \Theta = \{ \theta \in \mathbb{R}^{n_{p}} | H_{\theta} \theta \leq h_{\theta} \},
    \end{equation}
    with $H_{\theta}\in\mathbb{R}^{q_{\theta}\times n_{p}}$ and $h_{\theta}\in\mathbb{R}^{q_{\theta}}$.  
    The disturbance is bounded such that $\bar{M} W_j \in \mathbb{W}_x$, and $\bar{N} W_j\in \mathbb{W}_y $, $W_j \in \mathbb{W}^{p}$, where $\mathbb{W}_x$, $\mathbb{W}_y$ are known bounded polytopic sets. 
\end{assumption}
\begin{remark}\label{rem:uncertaintyass}
    This uncertainty characterization also captures the multi-step predictor characterization in~\cite{terzi2019learning,terzi2022robust} and contains impulse response models~\cite{tanaskovic2014adaptive} and orthonormal basis functions~\cite{heirung2017dual} as special cases.
    Assumption~\ref{ass:bound} can be directly derived from data within a set-membership framework, as shown in works like~\cite{lauricella2020set} and~\cite{terzi2022robust}. 
    These approaches result in $\mathbb{W}_x$, $\mathbb{W}_y$, and $\Theta$ in the shape of hyperboxes with $n_p=pn_x^2+n_xn_u \frac{p(p+1)}{2}$ and can even account for noisy output measurements.
    By directly identifying multi-step predictors we can verify smaller prediction error bounds, as demonstrated in recent works~\cite{terzi2022robust,terzi2019learning,lauricella2020set}.
    This reduces the conservatism of our approach, as discussed in Sec.~\ref{S:discussion}. 
\end{remark}
\begin{problem} \label{pr:controlProb}
    The control objective is to find a MPC formulation that stabilizes system~\eqref{eq:system} and ensures constraint satisfaction by utilizing the multi-step predictor with uncertainty bounds from Assumption~\ref{ass:bound}.
\end{problem}
\section{Multi-rate homothetic MPC}
In this section, we describe the proposed MPC formulation. 
Firstly, we describe the multi-rate tube MPC scheme in Sec.~\ref{SS:multirate}.
Afterward, we show how to construct a homothetic tube for the multi-step formulation in Sec.~\ref{SS:tube}. 
We finally introduce the overall formulation of the MPC problem and its theoretical analysis in Sec.~\ref{SS:formulation}.
\subsection{Conceptual multi-rate tube MPC for multi-step predictors} \label{SS:multirate}
Due to inconsistencies between the models, the multi-step model defined in~\eqref{eq:multistepmodel} cannot be directly applied in a standard MPC approach.
Therefore, non-standard MPC formulations, as been employed in~\cite{maddalena2021kpc,saccani2022multitrajectory} and~\cite{terzi2022robust}. In this work, we adopt the approach outlined in~\cite{terzi2022robust} by implementing a multi-rate MPC scheme.
We consider  a prediction horizon $N=N_p p$, which is split in $N_p$ intervals of length $p$.
Furthermore, the optimal control problem is only re-solved at any index $j$, i.e., every $p$ steps. 
Within the interval  $k\in\{jp,\ldots,(j+1)p-1\}$ the first inputs $U_j\in\mathbb{U}^p$ are applied. 
Thus, we essentially utilize the $p$-stead ahead predictor~\eqref{eq:multistepmodel} like a state space model with a slower sampling rate.

To solve the considered control problem, a tube-based MPC algorithm is considered, which predicts a sequence of sets $\{\pazocal{X}_{l|j}\}_{l\in\mathbb{N}_0^{N_{p}}}$ that provide an outer bound for the uncertain predicted state trajectory.
Specifically, these tubes are constructed  using a feedback policy $U_{l|j}(X)$ and the uncertainty bounds $\Theta$, $\mathbb{W}_x$ and $\mathbb{W}_y$. Enforcing state and input constraint satisfaction~\eqref{eq:constraints} over a horizon $N_p$ with such a tube formulation can be stated as:
\begin{subequations} \label{eq:tube}
\begin{align}
&\pazocal{X}_{0|j} \ni X_j \label{seq:parametrizationTube_currState}\\
&\pazocal{X}_{l+1|j} \ni \bar{A}(\theta) X+ \bar{B}(\theta) U_{l|j}(X) + \bar{M} W , \label{seq:parametrizationTube_stateEv}\\
& \ \ \ \ \forall \bar{M}W \in \mathbb{W}_x, \ \theta \in \Theta, \  X\in \pazocal{X}_{l|j}, \  l\in \mathbb{N}_0^{N_{p}-1} \nonumber
\\
&x \in \mathbb{X}, \  U_{l|j}(X)\in \mathbb{U}^{p},  \ \ \ \ \ \ \forall X\in \pazocal{X}_{l|j}, \ l\in \mathbb{N}_0^{N_{p}-1} \label{seq:parametrizationTube_constraints}\\
&\mathbb{X}^{p-1} \ni \bar{C} (\theta) X+\bar{D} (\theta) U_{l|j}(X) + \bar{N} W, \label{seq:parametrizationTube_output}\\
& \ \ \ \ \ \ \forall \bar{N}W \in \mathbb{W}_y, \ \theta \in \Theta, \  X\in \pazocal{X}_{l|j}, \  l\in \mathbb{N}_0^{N_{p}-1}.\nonumber
\end{align}
\end{subequations}
Specifically, Equs.~\eqref{seq:parametrizationTube_currState}-\eqref{seq:parametrizationTube_stateEv} ensure that uncertain state predictions lie in tube, while Equs.~\eqref{seq:parametrizationTube_constraints}-\eqref{seq:parametrizationTube_output} ensure that state and input constraints are satisfied at each time $k\in \mathbb{N}$.
At each time step $j$ (
i.e., every $p$ time instants), we solve the following finite horizon optimal control problem:
\begin{align} \label{eq:GeneralFHOCP}
    & \min_{\pazocal{X}_{\cdot|j}, U_{\cdot|j}(\cdot)} J_N(X_j,U_{\cdot|j}(\cdot)) \\
     \text{s.t.}&~\eqref{eq:tube}, \ \  \pazocal{X}_{{N_{p}}|j}\subseteq \pazocal{X}_{\mathrm{f}}, \nonumber
\end{align}
where $J_{N}$ is a cost function and $\pazocal{X}_{\mathrm{f}}$ a polytopic terminal set, which are specified later. 
Problem~\eqref{eq:GeneralFHOCP} provides a sequence of feedback policy that guarantee constraints satisfaction. At each time step $j$ (every $p$ steps), we implement the first component of the minimizer, denoted as $U(X_j)=U^\star_{0|j}$. 
However, it is not possible to solve Problem~\eqref{eq:GeneralFHOCP} without a tractable parametrization of both the feedback policy $U(X)$ and the state tube $\pazocal{X}_{\cdot|j}$.
%
\subsection{Homothetic tube parametrization} \label{SS:tube}
In the following, we introduce a finite-dimensional parametrization for Problem~\eqref{eq:GeneralFHOCP} by adapting the homothetic tube MPC in~\cite{lorenzen2019robust}.
Let us consider the following parametrization of the feedback policy $U(X)$:
\begin{equation} \label{eq:inputparametrization}
    U_{l|j}(X)=KX+V_{l|j},
\end{equation}
where $\bm{V}_{\cdot|j}=\{ V_{l|j} \}_{l\in \mathbb{N}_0^{{N_{p}}-1}} $, is  a sequence of decision variables and $K\in\mathbb{R}^{pn_u\times n_x}$ is an offline chosen stabilizing gain.
The feedback $K$ in combination with a Lyapunov function can be computed by solving a Linear Matrix Inequality (LMI), see Appendix~\ref{A:LMI} for details.
Let us now consider a given polytope $\pazocal{X}_0=\{ x\in\mathbb{R}^{n_x} | H_x x \leq \bm{1} \}$, where $H_x\in\mathbb{R}^{q_x\times n_x}$, and the MPC decision variables $\bm{z}_{\cdot|j}=\{ z_{{l|j}}\}_{l\in\mathbb{N}_0^{N_{p}}}$, $z_{l|j}\in\mathbb{R}^{n_x}$, and $\bm{\alpha}_{\cdot|j}=\{ \alpha_{l|j} \}_{l\in\mathbb{N}_0^{N_{p}}}$, $\alpha_{l|j} \in \mathbb{R}_{\geq 0}$. \\ The tube is parametrized as follows:
\begin{align}
    \pazocal{X}_{l|j}=&\{z_{l|j} \} \oplus \alpha_{l|j} \pazocal{X}_0 \nonumber \\
    =& \{x\in\mathbb{R}^{n_x} \ | \ H_x(x-z_{l|j})\leq \alpha_{l|j}\bm{1}\} \\
    = & \{z_{l|j}\} \oplus \alpha_{l|j} \ \text{co}\{ x^1,x^2,\dots x^{n_v} \} , \nonumber
\end{align}
where $x^v$, $v=1,\dots,n_v$ are the vertices of $\pazocal{X}_0$.\\
We define the ensuing quantities:
\begin{align}
    E^v_{l|j} =&\begin{bmatrix}
        \bar{A}_1X^v_{l|j}+\bar{B}_1U^v_{l|j},\dots , \bar{A}_{n_{p}}X^v_{l|j}+\bar{B}_{n_{p}}U^v_{l|j},\\
        \bar{C}_1X^v_{l|j}+\bar{D}_1U^v_{l|j},\dots ,\bar{C}_{n_{p}}X^v_{l|j}+\bar{D}_{n_{p}}U^v_{l|j}
    \end{bmatrix}, \label{seq:Ev} \\
     e_{l|j}^v = & \begin{bmatrix}
         \bar{A}_0X^v_{l|j} + \bar{B}_0U^v_{l|j} - z_{l+1|j}, \\
         \bar{C}_0X^v_{l|j}+\bar{D}_0U^v_{l|j}
     \end{bmatrix} ,\label{seq:ev}\\
     X^v_{l|j}= & z_{l|j}+\alpha_{l|j} x^v, \ \ \  U^v_{l|j}=V_{l|j}+KX^v_{l|j}, \label{seq:XvUv}\\
     [\bar{f}]_i = & \max_{x\in\pazocal{X}_0}[F]_i x, \  \ \ [\bar{g}]_j =  \max_{x\in\pazocal{X}_0}[G_{p}K]_j x,\label{seq:maximizationStateConstr}\\
     & \ \ \forall i \in \mathbb{N}_0^{c_{\mathrm{x}}}, \ \ \ \ \forall j \in \mathbb{N}_0^{c_{\mathrm{u}}p}, \nonumber \\
     [\bar{w}_x]_i =& \max_{w_x\in\mathbb{W}_x} [H_x]_i w_x,  \ \ \  [\bar{w}_y]_j = \max_{w_y\in\mathbb{W}_y} [F_{p}]_j w_y \\
     & \ \ \forall i \in \mathbb{N}_0^{q_x}, \ \ \ \ \forall j \in \mathbb{N}_0^{c_{\mathrm{x}}(p-1)}, \nonumber\\
      H_{p}=&\text{blkdiag}\{ H_x, F_{p} \}, \ \ \  \bar{w}=[w_x^{\top}, \ w_y^{\top}]^{\top} \\
      \Phi_{l+1|j} = &[\alpha_{l+1|j}\bm{1}^{\top}, \bm{1}^{\top}]^{\top}\in\mathbb{R}^{q_x+c_{\mathrm{x}}(p-1)},
\end{align}
where $F_{p}=\text{blkdiag}\{F, \dots, F \} \in\mathbb{R}^{c_{\mathrm{x}} (p-1)\times n_x (p-1) }$ and $G_{p}=\text{blkdiag}\{G, \dots, G \} \in\mathbb{R}^{c_{\mathrm{u}} p\times n_up }$. 
The following proposition introduces a non-conservative reformulation of the tube constraints~\eqref{eq:tube}, which accounts for the parametric uncertainty $\theta\in\Theta$.
\begin{proposition}{(Homothetic tube parametrization)}
\label{prop:tube}
    \eqref{seq:parametrizationTube_currState}-\eqref{seq:parametrizationTube_output} are satisfied if and only if there exists \quad\quad\quad\quad\quad\quad \vspace{0mm} $\Lambda^v_{l|j}\in\mathbb{R}_{\geq0}^{q_x+c_{\mathrm{x}} (p-1)\times q_{\theta}}$ such that \\
\begin{subequations} \label{eq:tube_{p}arametrization}
    \begin{align}
    &F z_{l|j} +\alpha_{l|j}\bar{f}\leq \bm{1}, \label{seq:constraints_state} \\
    &G_{p}K z_{l|j} + G_{p} V_{l|j} + \alpha_{l|j}\bar{g}\leq \bm{1}, \label{seq:constraints_input} \\
    &H_x(x_k-z_{0|j})\leq \alpha_{0|j}\bm{1}, \label{seq:currentState} \\
    & \Lambda_{l|j}^v h_{\theta}+H_{p} e_{l|j}^v \leq  \Phi_{l+1|j} -\bar{w}, \label{seq:prediction1}\\ 
    & H_{p} E_{l|j}^v=\Lambda_{l|j}^v H_{\theta},\label{seq:prediction2}
\end{align}
\end{subequations}
hold for all $v\in\mathbb{N}_1^{n_v}$, $l\in\mathbb{N}_0^{N_p-1}$.
\end{proposition}
\begin{proof}
The proof follows similar steps to the one provided in~\cite[Prop. 9]{lorenzen2019robust}. 
The inequality~\eqref{seq:parametrizationTube_constraints} is equivalent to  
\begin{align*}
    &F z_{l|j} +\alpha_{l|j}Fx \leq \bm{1} \ \ \ \ &&\forall x \in\pazocal{X}_0\\
    &G_{p}K z_{l|j} + G_{p} V_{l|j} + \alpha_{l|j}G_{p}Kx\leq \bm{1}&&\forall x \in\pazocal{X}_0
\end{align*}
which are equivalent to~\eqref{seq:constraints_state},~\eqref{seq:constraints_input} using the maximization~\eqref{seq:maximizationStateConstr}. Inequality~\eqref{seq:parametrizationTube_currState} is equivalent to~\eqref{seq:currentState}. 
Using the vertices~\eqref{seq:XvUv}, Equations~\eqref{seq:parametrizationTube_stateEv} and ~\eqref{seq:parametrizationTube_output} are equivalent to
\begin{align*}
    &\left\{\begin{array}{l}\pazocal{X}_{l+1|j}\ominus\mathbb{W}_{\mathrm{x}} \ni \bar{A}(\theta)X^v_{l|j}+\bar{B}(\theta) U^v_{l|j}  \\
    \mathbb{X}^{p-1} \ominus\mathbb{W}_{\mathrm{y}}\ni \bar{C}(\theta)X_{l|j}^v + \bar{D}(\theta) U^v_{l|j} \ \ \  \end{array}\right\} \begin{matrix}
        \forall v \in \mathbb{N}_1^{n_v}, \\ \theta \in \Theta
    \end{matrix}\\
    &\Leftrightarrow \left\{\begin{array}{l} H_{\mathrm{x}} ( \bar{A}(\theta) X_{l|j}^v+\bar{B}(\theta) U^v_{l | k}\\
    -z_{l+1 | k} ) \leq \alpha_{l+1|j} \mathbf{1}-\bar{w}_{\mathrm{x}}, \\
    F_{p}(\bar{C}(\theta)X_{l|j}^v+\bar{D}(\theta) U^v_{l | k})  \leq \mathbf{1}-\bar{w}_{\mathrm{y}}, 
    \end{array}\right\} \begin{matrix}
        \forall v \in \mathbb{N}_1^{n_v}, \\ \theta \in \Theta
    \end{matrix} \\ 
    &\Leftrightarrow  \max_{\theta\in\Theta}H_{p}E^v_{l|j}\theta+H_{p}e^v_{l|j}\leq \Phi_{l+1|j}-\bar{w} ,~ \forall v \in \mathbb{N}_1^{n_v}
     \\ 
    &\Leftrightarrow\left\{\begin{array}{l}\Lambda_{l|j}^v h_{\theta}+H_{p} e_{l|j}^v \leq  \Phi_{l+1|j} -\bar{w} \\ 
    H_{p} E_{l|j}^v=\Lambda_{l|j}^v H_{\theta} \\ 
    \Lambda_{l|j}^v \in \mathbb{R}_{\geq 0}^{q_x+c_{\mathrm{x}} (p-1) \times q_{\theta}}\end{array}\right\} \quad \forall v \in \mathbb{N}_1^{n_v},
\end{align*}
which corresponds to~\eqref{seq:prediction1}--\eqref{seq:prediction2}, where the last equivalence utilized strong duality of linear programming exploiting half-space representation of $\Theta$~\eqref{eq:modelUncer}.
\end{proof}
\subsection{MPC algorithm and theoretical analysis} \label{SS:formulation}
We are now in position to define the MPC problem solved at each time step $j$.
For simplicity, we consider a nominal cost using a prediction of $\hat{X}_{\cdot|j}$ and $\hat{Y}_{\cdot|j}$ based on some nominal parameter $\hat{\theta}\in\Theta$, i.e.,
\begin{subequations} \label{eq:nominalPred}
\begin{align}
    \hat{X}_{0|j} =& X_j \label{eq:FHOCP_noinal_start}, &&\\
    \hat{X}_{l+1|j} =&\bar{A}(\hat{\theta}) \hat{X}_{l|j} +\bar{B}(\hat{\theta}) \hat{U}_{l|j}, &&\\
    \hat{Y}_{l|j} =& \bar{C}(\hat{\theta}) \hat{X}_{l|j} +\bar{D}(\hat{\theta}) \hat{U}_{l|j}, &&\\
    \hat{U}_{l|j} =& K \hat{X}_{l|j} + V_{l|j}, && \forall l\in\mathbb{N}_0^{N_p-1}\label{eq:FHOCP_noinal_end} .
\end{align}
\end{subequations}
The MPC cost function is given by:
\begin{multline} \label{eq:cost}
J_{N}(\hat{\bm{X}}_{\cdot|j},\hat{\bm{Y}}_{\cdot|j},\hat{\bm{U}}_{\cdot|j}) = \sum_{l=0}^{N_{p}-1} ( \| \hat{Y}_{l|j} \|_{Q_{p}}^2 +\\ + \| \hat{X}_{l|j} \|_Q^2 + \| \hat{U}_{l|j} \|_{R_{p}}^2  ) + \| \hat{X}_{N_{p}|j} \|_{P}^2.
\end{multline}
$Q\in\mathbb{R}^{n_x\times n_x}$ and $R\in\mathbb{R}^{n_u\times n_u}$ are user-chosen positive definite matrices that can be tuned to achieve desired performances~\cite{kouvaritakis2016model} and ${Q_{p}}=\text{blkdiag}\{ Q,\dots,Q\}\in\mathbb{R}^{n_x(p-1)\times n_x(p-1)}$, ${R_{p}}=\text{blkdiag}\{ R,\dots,R\}\in\mathbb{R}^{n_u p\times n_up}$. The matrix $P\in\mathbb{R}^{n_x\times n_x}$ corresponds to a Lyapunov function (see Asm.~\ref{ass:terminal}).
To simplify the notation, let us define the following set of decision variables $\bm{o}_{\cdot|j} = \{ \bm{z}_{\cdot|j}, \bm{\alpha}_{\cdot|j}, \bm{V}_{\cdot|j}, \bm{\Lambda}_{\cdot|j}, \bm{\hat{X}}_{\cdot|j}, \bm{\hat{U}}_{\cdot|j}, \bm{\hat{Y}}_{\cdot|j} \}$, with $\bm{\Lambda}_{\cdot|j}=\{ \Lambda^v_{l|j} \}_{v\in\mathbb{N}_1^{n_v}, l\in\mathbb{N}_0^{N_{p}}}$.
Thus, the proposed MPC formulation is given by the following linearly constrained quadratic program (QP):
\begin{subequations} \label{eq:FHOCP}
\begin{align}
        \bm{o}^\star_{\cdot|j}= &\argmin_{\bm{o}_{\cdot|j}} J_{N}(\hat{\bm{X}}_{\cdot|j},\hat{\bm{Y}}_{\cdot|j},\hat{\bm{U}}_{\cdot|j}) \\
        \text{s.t.  }& 
         z_{{N_{p}}|j}\oplus\alpha_{{N_{p}}|j}\pazocal{X}_0 \subseteq \pazocal{X}_{\mathrm{f}} \label{eq:FHOCP_terminal} \\
        &\eqref{seq:Ev}, \eqref{seq:ev}, \eqref{seq:XvUv}, \eqref{eq:tube_{p}arametrization},\eqref{eq:nominalPred}.\label{eq:FHOCP_end}
\end{align}
\end{subequations}
We solve Problem~\eqref{eq:FHOCP} at each time step and $j$ (every $p$ time instants) and apply the input $U_j=V_{0|j}^\star+KX_j$.
\begin{remark}
In many cases, the parameter set $\Theta$ is structured. For example, in~\cite{lauricella2020set}, each row of $\bar{A},\ \bar{B}, \ \bar{C}, \ \bar{D}$ is parameterized with independent parameters. This structured arrangement leads to a sparsity pattern in the dual variables $\Lambda^v_{l|j}$, which significantly reduces the number of decision variables. For further details, refer to the accompanying code in Sec.~\ref{S:num_example}.
\end{remark}
We consider the following assumption for the terminal cost $P$ and terminal set $\pazocal{X}_{\mathrm{f}}$.
\begin{assumption}[Terminal set] \label{ass:terminal}
    It holds $\pazocal{X}_{\mathrm{f}}=\eta\pazocal{X}_0$ with $\eta>0$.
    The set $\pazocal{X}_{\mathrm{f}}$ is a Robust Positevely Invariant (RPI) set for system~\eqref{eq:multistepmodel} with the control law $U=KX\in\mathbb{U}^{p}$ and ensures constraint satisfaction, i.e.,
    $\bar{A}_{\mathrm{cl}}(\theta) \pazocal{X}_{\mathrm{f}}\oplus \mathbb{W}_{\mathrm{x}}\subseteq \pazocal{X}_{\mathrm{f}}$, $K\pazocal{X}_{\mathrm{f}}\subseteq \mathbb{U}^{p}$ and $\bar{C}_{\mathrm{cl}}(\theta)\pazocal{X}_{\mathrm{f}}\oplus\mathbb{W}_{\mathrm{y}}\subseteq\mathbb{X}^{p-1}$, $\forall \theta\in\Theta$, where $\bar{A}_{\mathrm{cl}}(\theta)=\bar{A}(\theta)+\bar{B}(\theta)K$ and $\bar{C}_{\mathrm{cl}}(\theta)=\bar{C}(\theta)+\bar{D}(\theta)K$.
\end{assumption}
Furthermore, the matrix $P\in\mathbb{R}^{n_x\times n_x}$ in~\eqref{eq:cost} satisfy
\begin{multline} \label{eq:Riccati}
    \bar{A}_{\mathrm{cl}}(\theta)^{\top} P \bar{A}_{\mathrm{cl}}(\theta) + Q +\bar{C}_{\mathrm{cl}}(\theta)^{\top} Q_{p} \bar{C}_{\mathrm{cl}}(\theta)+ \\
    + K^{\top}R_pK \preceq P \ \ \ \  \forall \theta \in \Theta.
\end{multline}
Inequality~\eqref{eq:Riccati} can be reformulated as an LMI.
In Appendix~\ref{A:LMI} we propose a reformulation of the LMI exploiting the structure of $\Theta$ to avoid enumerating all vertices of the high-dimensional set $\Theta$.

The following theorem establishes recursive feasibility and practical asymptotic stability of the proposed MPC scheme.
\begin{theorem}
    Let Assumptions~\ref{ass:bound} and~\ref{ass:terminal} hold, and suppose that at time $j=0$ Problem~\eqref{eq:FHOCP} is feasible. Then, the MPC scheme~\eqref{eq:FHOCP} is feasible for all time instant $j>0$, the resulting closed-loop system robustly satisfies constraints~\eqref{eq:constraints} for all $k\in\mathbb{N}$ and $X_j=0$ is practically asymptotically stable~\cite[Def.~2.2]{grune2014asymptotic}.
\end{theorem}
\begin{proof}
\textit{Recursive feasibility: }
    The constraints~\eqref{eq:FHOCP_noinal_start}--\eqref{eq:FHOCP_noinal_end} are always feasible and
      due to Proposition~\ref{prop:tube} the constraints~\eqref{eq:FHOCP_terminal}--\eqref{eq:FHOCP_end}  are equivalent to Problem~\eqref{eq:GeneralFHOCP}, i.e., it suffices to study recursive feasibility of Problem~\eqref{eq:GeneralFHOCP} under the given parameterization.
    We conduct the proof by induction, i.e., assuming Problem~\eqref{eq:GeneralFHOCP} is feasible at some time $j$, we construct a feasible candidate solution at time $j+1$. At time $j$, the optimal solution is denoted by $\pazocal{X}^\star_{l|j}$ with feedback $U_{l|j}^\star(X)$, $l\in\mathbb{N}_0^{N_p}$.  
    The feasible candidate sequence at time $j+1$ is given by $\pazocal{X}_{l|j+1}=\pazocal{X}_{l+1|j}^\star$, $U_{l|j+1}(X)=U_{l+1|j}^\star(X)$, $l\in\mathbb{N}_0^{N_p-1}$ and $\pazocal{X}_{N_p|j+1}=\pazocal{X}_{\mathrm{f}}$, 
    $U_{N_p|j+1}(X)=KX.$    
  Equ.~\eqref{seq:parametrizationTube_currState} holds since $X_{j+1}\in\pazocal{X}_{1|j}^\star=\pazocal{X}_{0|j+1}$ using the assumed disturbance bound (Asm.~\ref{ass:bound}) and that $\pazocal{X}_{1|j}^\star$ satisfies~\eqref{seq:parametrizationTube_stateEv} for $l=0$ and $X_j\in\pazocal{X}_{0|j}^\star$. 
 Conditions~\eqref{seq:parametrizationTube_stateEv}, \eqref{seq:parametrizationTube_constraints}, \eqref{seq:parametrizationTube_output} hold directly for $l\in\mathbb{N}_0^{N_p-2}$ due to the shifting. 
 Assumption~\ref{ass:terminal} implies satisfaction of \eqref{seq:parametrizationTube_stateEv}--\eqref{seq:parametrizationTube_output}
 for $l=N_p-1$. 
The terminal set constraint is satisfied by definition. \\
\textit{Constraint satisfaction:} 
Feasibility of Problem~\eqref{eq:tube} with~\eqref{seq:parametrizationTube_currState},~\eqref{seq:parametrizationTube_constraints},~\eqref{seq:parametrizationTube_output} for $l=0$ ensures that the resulting state and input satisfy $x_k\in\mathbb{X}$, $u_k\in\mathbb{U}$ for $k=jp,\dots,(j+1)p-1$, i.e., the closed-loop system satisfies the constraints~\eqref{eq:constraints}. \\
    \textit{Stability:}
    The stability proof follows standard arguments from robust MPC~\cite{kouvaritakis2016model}. 
    Let us denote the minimum in~\eqref{eq:FHOCP} by $J_N^\star(X_j)$ for a given state $X_j$. 
    In case of exact predictions, i.e., $X_{j+1}=\hat{X}_{1|j}^\star$, the candidate solution in combination with the terminal cost~\eqref{eq:Riccati} implies
    \begin{multline*}
    J_N^\star(\hat{X}_{1|j}^\star)-J_N^\star(X_j)\\
    \leq-\left(\|\hat{Y}_{0|j}^\star\|_{Q_{p}}^2+\|\hat{X}_{0|j}^\star\|_Q^2+\|\hat{U}_{0|j}^\star\|_{R_{p}}^2\right)\leq -\|X_j\|_Q^2.   
    \end{multline*}
    The value function of the QP~\eqref{eq:FHOCP} is  piecewise quadratic and hence uniformly continuous 
    , i.e., there exists a function $\gamma\in\mathcal{K}$ such that $|J_N^\star(x)-J_N^\star(\tilde{x})|\leq \gamma(\|x-\tilde{x}\|)$ for any two states $x,\tilde{x}$ for which~\eqref{eq:FHOCP} is feasible. 
    Given recursive feasibility, it holds that
    \begin{align*}
    &J_N^\star(X_{j+1})\leq 
    J_N^\star(X_{1|j}^\star)+\gamma(\|X_{j+1}-\hat{X}_{1|j}^\star\|)\\
    \leq& J_N^\star(X_j)-\|X_j\|_Q^2+\gamma(\|X_{j+1}-\hat{X}_{1|j}^\star\|).
    \end{align*}
    The prediction error $X_{j+1}-X_{1|j}^\star=(\bar{A}(\theta^*)-\bar{A}(\hat{\theta}))X_j+(\bar{B}(\theta^*)-\bar{B}(\hat{\theta}))U_j+\bar{M} W_j$ is bounded since $\mathbb{X},\Theta,\mathbb{U},\mathbb{W}$ are bounded, i.e., there exist a constant $\epsilon>0$ s.t. 
    \begin{align*}
    &J_N^\star(X_{j+1})\leq J_N^\star(X_j)-\|X_j\|_Q^2+\gamma(\epsilon).
    \end{align*}
    Furthermore, $J_N^\star(X)$ piece-wise quadratic ensures that it is lower and upper bound by a $\mathcal{K}_\infty$ function of $\|X\|$.
    Thus,  $J_N^\star$ is a practical Lyapunov function~\cite[Def.~2.3]{grune2014asymptotic}, which implies practical asymptotic stability of $X_j=0$~\cite[Thm.~2.4]{grune2014asymptotic}. 
\end{proof}

\section{Discussion}\label{S:discussion}
In this section, we compare and contrast the proposed multi-rate homothetic tube MPC with state-of-the-art approaches. 
Specifically, we focus on the multi-rate MPC presented in \cite{terzi2022robust} and the homothetic tube MPC introduced in \cite{lorenzen2019robust}. 
In the numerical comparison (Section~\ref{S:num_example}), we additional compare~\cite{Koehler2019Adaptive} which provides an alternative homothetic tube formulation.
Firstly, in the special case $p=1$, Problem~\eqref{eq:FHOCP} reduces to the homothetic tube MPC described in~\cite{lorenzen2019robust}. The key differentiating factor of our approach is the incorporation of a multi-step predictor~\eqref{eq:multistepmodel} in combination with a multi-rate formulation.
Secondly, the considered multi-rate formulation and the uncertainty description utilizing multi-step predictors (Asm.~\ref{ass:bound}) is based on~\cite{terzi2022robust}.
A fundamental distinction lies in the prediction of the tube $\pazocal{X}$.
The homothetic tube formulation non-conservatively accounts for the parametric uncertainty $\theta\in\Theta$ (Prop.~\ref{prop:tube}), while \cite{terzi2022robust} utilizes a constant $\tau_p\in\mathbb{R}^{n_x}$ 
\begin{align}
\label{eq:tau}    
\tau_{p}\geq |\bar{A}(\hat{\theta})X_j+\bar{B}(\hat{\theta})U_j+\bar{M}W_j-X_{j+1}|, 
\end{align}
which upper bounds the prediction error for all possible state and input pairs $X_j,U_j,W_j$\footnote{%
This bound can be computed using Asm.~\ref{ass:bound} and compact constraints $X_j\in\mathbb{X}$, $U_j\in\mathbb{U}^{p}$.
In~\cite{lauricella2020set,terzi2019learning,terzi2022robust}, a less conservative bound is computed by using the pairs $(X_j,U_j)$ observed in the offline data. However, resulting robustness guarantees may fail if closed-loop operation differs significantly from the observed data.}.
In particular, the multi-rate MPC~\cite{terzi2022robust} utilizes a (standard/rigid) tube MPC~\cite{mayne2005robust} with $\pazocal{X}_{l|j}=\{z_{l|j}\}\oplus\pazocal{X}_0$, where $\pazocal{X}_0$ is an RPI set w.r.t. the disturbance bound $\tau_{p}$~\eqref{eq:tau}.
\subsubsection*{Conservatism}
Compared to \cite{terzi2022robust}, our approach is inherently less conservative (given the same tube parametrization $\pazocal{X}_0$ and feedback $K$) due to the utilization of a homothetic tube, which provides an exact $p$-step reachable set (Prop.~\ref{prop:tube}). In contrast, $\tau_{p}$~\eqref{eq:tau} used in~\cite{terzi2022robust} represents a conservative over-approximation considering the worst-case model mismatch in the full range operation. 
In the next section, we demonstrate this reduced conservatism numerically, using the example provided in~\cite{terzi2022robust}. The primary advantage over~\cite{lorenzen2019robust} and other homothetic MPC approaches in the literature lies in the utilization of multi-step predictors.
Multi-step predictors allow for simple bounds on the prediction error and enable a more straightforward reformulation to address parametric uncertainty. To exploit these benefits, Proposition~\ref{prop:tube} enables an exact reformulation over the horizon $p$ of the multi-step predictor. Conversely, addressing parametric uncertainty in state space models usually necessitates sequential propagation~\cite{kouvaritakis2016model,lorenzen2019robust,Lu2019RAMPC,Koehler2019Adaptive,parsi2022scalable,langson2004robust}. 
Specifically, in~\cite{lorenzen2019robust}, each time step requires a tube propagation that results in conservatism due to the fixed parametrization of $\pazocal{X}_0$.
In contrast, our multi-rate formulation requires only $1/p$ such over-approximations.
As an extreme case, the tube propagation in~\cite{lorenzen2019robust} becomes unbounded in case there exist no common Lyapunov function for $\bar{A}_{\mathrm{cl}}(\theta)$, $\theta\in\Theta$.
However, this limitation becomes negligible for multi-step predictors as $\|\bar{A}_{\mathrm{cl}}\|\rightarrow 0$ when $p\rightarrow \infty$ for any stable system. 
\subsubsection*{Computational complexity}
In~\cite{terzi2022robust}, the computational complexity remains comparable to a nominal MPC, necessitating only the inclusion of extra linear inequality constraints for the initial state related to the RPI set.
To streamline our discussion, we suppose each entry in $\bar{A}, \ \bar{B}, \ \bar{C}$ and $\bar{D}$ is characterized by a different parameter and $\Theta$ is a box constraint, see Rem.~\ref{rem:uncertaintyass}.
For the proposed multi-step predictor, we set the tube $\pazocal{X}_0$ as a low-complexity polytope~\cite[Sec.~5.4]{kouvaritakis2016model}, with $n_v=2^{n_x}$ vertices and $2n_x$ half-space constraints.
This selection is always feasible for open-loop stable systems, provided that $p$ is suitably large.
The number of constraints required to satisfy~\eqref{seq:parametrizationTube_stateEv} scales linearly with both the vertices $n_v=2^{n_x}$ and $n_x^2+n_xn_up$ half-space constraints of $\Theta$. 
In contrast, the homothetic tube MPC detailed in \cite{lorenzen2019robust}, tailored for state space models, involves $n_v$ vertices within $\pazocal{X}_0$ and $n_x^2+n_un_x$ half-space constraints on $\Theta$. Notably, the number of vertices of the polytopic tube $\pazocal{X}_0$ tend to scale exponentially with respect to $n_x$ and this can dominate the complexity for systems with a large number of states. 
As a result, our proposed approach shows potential for improved computational efficiency compared to \cite{lorenzen2019robust}, particularly when dealing with a large number of vertices in the polytopic tube $\pazocal{X}_0$, which can be computationally demanding, see Tab.~\ref{tab:example}.
\begin{remark}
    In the special case $N=p$, recursive feasibility with multi-step predictors can also be ensured with a simple constraint tightening without tube propagation by utilizing an adaptive horizon in the MPC~\cite{saccani2022multitrajectory}.
\end{remark}


\section{Numerical example} \label{S:num_example}
\begin{figure}
	\centering
    \includegraphics[width=.95\columnwidth]{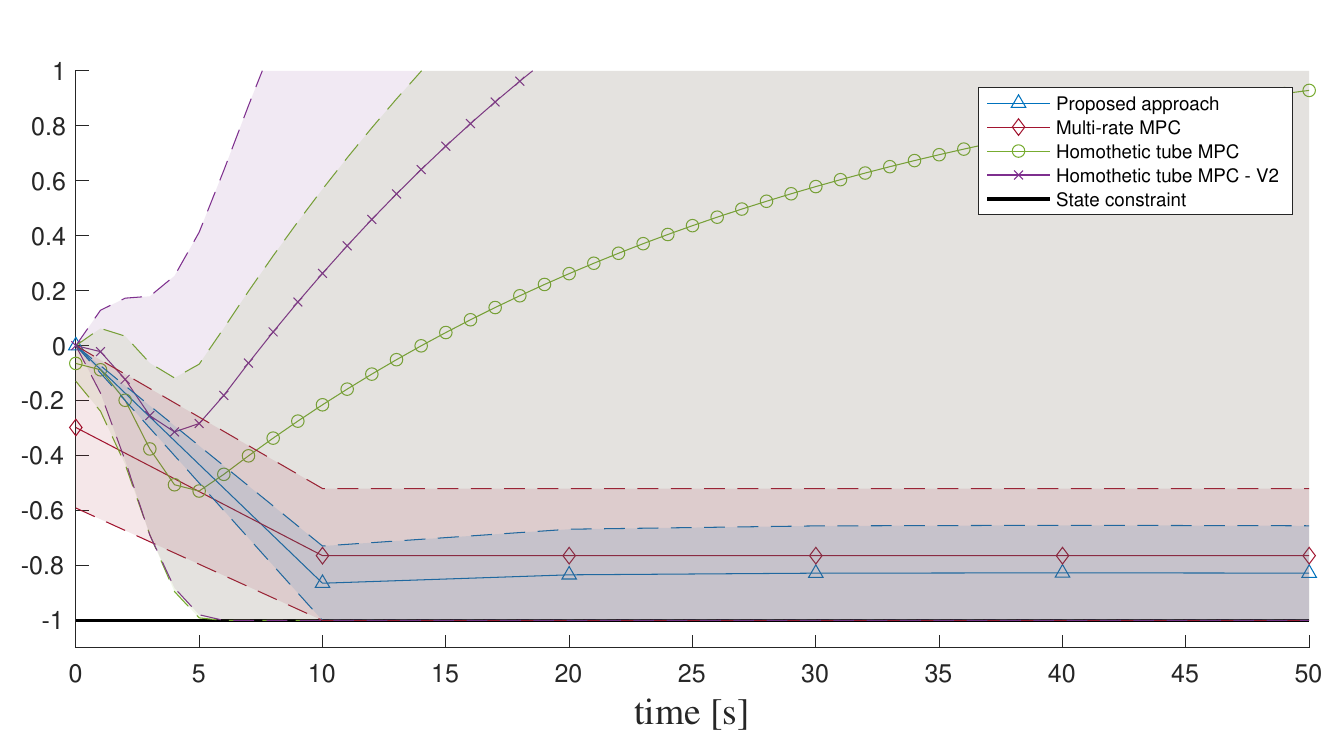}
	\caption{ Comparison of the open-loop solution of the predicted tube for the proposed approach (blue), the multi-rate MPC proposed in~\cite{terzi2022robust} (red) and the homothetic tube MPC~\cite{lorenzen2019robust} (green) and homothetic tube MPC - V2~\cite{Koehler2019Adaptive} (purple).
    The predicted sets are shaded with a dashed boundary and the nominal state prediction $z$ is solid with markers. The state constraint is highlighted in black.
    }
 \label{fig:output_tube}
 \vspace{-6mm}
\end{figure}
We highlight the benefits of our MPC framework over that of~\cite{terzi2022robust} by using the same numerical example.
The considered system is derived by discretizing the following continuous-time transfer function with a sampling time of $T_s = 0.1$:
\begin{equation}\label{eq:TF}
G(s) = \frac{160}{(s + 10)\left(s^2 + 1.6s + 16\right)}.
\end{equation}
Subsequently, the model is transformed into the discrete-time state space model~\eqref{eq:system} and the multi-step predictors~\eqref{eq:multistepmodel} identified by following the approach in~\cite{lauricella2020set}, with $p=10$.
Further details regarding the identification and uncertainty quantification can be found in Appendix~\ref{A:example} and the implementation of the proposed approach at: \url{https://github.com/DecodEPFL/HomotheticMPCmultistep}.
\subsubsection*{Comparison to state-of-the-art MPC formulations}
To demonstrate the reduced conservatism of the proposed approach we consider the objective of minimizing the value of the output of the transfer function while robustly ensuring a lower bound on the same value. 
This is formulated using $\mathbb{X}=\{x\in[-10,10]^3\ x_3\geq -1\}$, $\mathbb{U}=[-10,10]$, and $J_N$ minimizes $c^\top \hat{X}$, with $c=[0\ 0\ 1]^{\top}$. 
The proposed approach, along with~\cite{lorenzen2019robust} and~\cite{Koehler2019Adaptive}, utilizes as $\pazocal{X}_0$ a low-complexity polytope~\cite[Sec.~5.4]{kouvaritakis2016model}, while~\cite{terzi2022robust} employs a minimum RPI set. 
For ease of comparison, we exclusively examine the open-loop solution starting from an initial condition $x_0 = 0$ with a prediction horizon of $N_p = 5$, and we omit terminal constraints to ensure a fair comparison independent of the different terminal sets.
Fig.~\ref{fig:output_tube} illustrates the open-loop solutions of the proposed approach (blue), the multi-rate MPC introduced in~\cite{terzi2022robust} (red), and the homothetic tubes~\cite{lorenzen2019robust} (green) and~\cite{Koehler2019Adaptive} (purple), while Tab.~\ref{tab:example} provides details regarding the computational complexity of the QP for the different MPC formulations. 
The proposed approach results in a significantly smaller tube $\pazocal{X}$, enabling operation closer to the constraints. 
Compared to~\cite{terzi2022robust}, the approach is less conservative but results in increased computational complexity.
The approaches presented in~\cite{lorenzen2019robust,Koehler2019Adaptive} have a comparable computational complexity with a similarly simple tube $\pazocal{X}_0$; however, this results in significant conservatism, even an unbounded growth with the considered low-complexity polytope.  
This issue can be reduced using a more complex tube parametrization $\pazocal{X}_0$ e.g., an RPI set, but this would increase the number of vertices up to $10^3$ and number of decision variables/inequality constraints in~\cite{lorenzen2019robust} by a factor of $10^3$.
While alternative tube propagations in~\cite{Koehler2019Adaptive,Lu2019RAMPC} do not suffer from this limitation, they instead require enumeration of the vertices of $\Theta$, which limits scalability\footnote{%
In the considered example, the multi-step predictor has a parameter vector $\theta\in \mathbb{R}^{39}$, resulting in over $10^{10}$ vertices, which is not tractable.}.
Notably, both the proposed approach and the multi-rate MPC from~\cite{terzi2022robust} solve the optimization problem at each time step $j$, which corresponds to every $p=10$ time instants, while in~\cite{lorenzen2019robust,Koehler2019Adaptive} the optimization problem is solved at each time instant $k$. 
In conclusion, as shown in Fig.~\ref{fig:output_tube}, the proposed approach results in a tighter and thus less conservative tube compared to the other approaches, with moderate complexity (see Tab.~\ref{tab:example}).
\begin{table}
\vspace{2mm}
\begin{center}
\begin{tabular}{ | c | c | c | c| c| } 
  \hline
 Approach& \#opt. var. & \#ineq. con. & \#eq. con. \\ 
  \hline
  Multi-rate MPC~\cite{terzi2022robust} & $2.5\cdot10^2$ & $9.5\cdot10^2$ & $1.5\cdot10^2$  \\ 
  \hline
  Homothetic tube~\cite{lorenzen2019robust} & $7.1\cdot 10^4$ & $7.2\cdot 10^4$&$3.4\cdot 10^4$\\ 
  \hline
  Homothetic tube - V2~\cite{Koehler2019Adaptive} & $3.0\cdot 10^2$ &$1.2\cdot 10^6$ & $2.0\cdot 10^2$  \\ 
    \hline
  Proposed~\eqref{eq:FHOCP} & $1.2\cdot 10^5$ &$1.2\cdot 10^5$ & $6.1\cdot 10^5$  \\ 
  \hline

\end{tabular} \caption{Comparison of computational complexity of the QPs.}
\label{tab:example}
\end{center}
\vspace{-13mm}
\end{table}
\section{Conclusion}
We have introduced a homothetic tube MPC with multi-step predictors that combines the multi-rate MPC formulation utilizing multi-step predictors~\cite{terzi2022robust} with the reduced conservatism of homothetic tubes~\cite{lorenzen2019robust}. The benefits of this approach have been showcased through an in-depth comparison with the aforementioned methodologies and a numerical example.
The proposed approach paves the way for future research that explores its natural extensions, such as incorporating adaptive model updates~\cite{lorenzen2019robust,Lu2019RAMPC,Koehler2019Adaptive,tanaskovic2014adaptive} and addressing scalability limitations associated with vertex enumeration by utilizing ellipsoids~\cite{parsi2022scalable}.


\bibliographystyle{IEEEtran}
\bibliography{main}             

\appendix
\subsection{LMI derivation} \label{A:LMI}
To compute the control law $K$ and matrix $P$ satisfying~\eqref{eq:Riccati}, let us consider the following equation
\begin{multline}
\label{eq:terminal_cost}
    \bar{A}_{\mathrm{cl}}(\theta)^{\top} P \bar{A}_{\mathrm{cl}}(\theta) + Q +\bar{C}_{\mathrm{cl}}(\theta)^{\top} Q_p \bar{C}_{\mathrm{cl}}(\theta)+ \\
    + K^{\top}R_pK \preceq P \ \ \ \  \forall \theta \in \Theta.
\end{multline}
Suppose each entry in $\bar{A},\bar{B},\bar{C},\bar{D}$ corresponds to an uncertain parameter and $\Theta$ is a box, see Remark~\ref{rem:uncertaintyass}. 
Standard robust control designs for $P,K$ result in LMI conditions that scale with the number of vertices $\Theta$, which can be very large. 
In the following, we derive LMI conditions that avoid enumerating all vertices of the high-dimensional set $\Theta$. 
For ease of notation, we drop the explicit dependence on $\theta$ in the derivation and all conditions are imposed $\forall \theta\in\Theta$. 
Denote $\bar{C}_{\mathrm{cl}}=[\bar{C}_{\mathrm{cl},1}^\top,\dots \bar{C}_{\mathrm{cl},p-1}^\top]^\top$, with $\bar{C}_{\mathrm{cl},j}$ corresponding to prediction of $x$ $j$-time steps in the future.
It holds
\begin{align*}
\bar{C}_{\mathrm{cl}}^{\top} Q_p \bar{C}_{\mathrm{cl}}=\sum_{j=1}^{p-1} \bar{C}_{\mathrm{cl},j}^\top Q \bar{C}_{\mathrm{cl},j}.
\end{align*}
Condition~\eqref{eq:terminal_cost} is equivalent to existence of $\bar{Q}_j$ s.t.
\begin{align*}
&\bar{A}_{\mathrm{cl}}^{\top} P \bar{A}_{\mathrm{cl}} + Q+\sum_{j=1}^{p-1} \bar{Q}_j+  K^{\top}R_pK \preceq P    \\
&\bar{C}_{\mathrm{cl},j}^\top Q \bar{C}_{\mathrm{cl},j}\preceq \bar{Q}_j,~j\in\mathbb{N}_1^{p-1}.
\end{align*}
We assume that $Q,P$ are diagonal, i.e., $Q=\mathrm{diag}(q_1,\dots,q_{n_x})$ and $P=\mathrm{diag}(p_1,\dots, p_{n_x})$. 
It holds
\begin{align*}
\bar{A}_{\mathrm{cl}}^{\top} P \bar{A}_{\mathrm{cl}}&=\sum_{i=1}^{n_x}[\bar{A}_{\mathrm{cl}}]_i^\top p_i [\bar{A}_{\mathrm{cl}}]_i\\
\bar{C}_{\mathrm{cl},j}^{\top} Q \bar{C}_{\mathrm{cl},j}&=\sum_{i=1}^{n_x}[\bar{C}_{\mathrm{cl},{j}}]_i^\top q_i [\bar{C}_{\mathrm{cl},{j}}]_i,
\end{align*}
where $[\bar{C}_{\mathrm{cl},{j}}]_i$ corresponds to the $j$-step predictor for the $i$-th state component. 
Thus, condition~\eqref{eq:terminal_cost} is equivalent to existence of $\bar{Q}_{i,j}$, $\bar{P}_i$ satisfying 
\begin{align*}
&\sum_{i=1}^{n_x} \bar{P}_i+ Q+\sum_{j=1}^{p-1}\sum_{i=1}^{n_x} \bar{Q}_{i,j}+  K^{\top}R_pK \preceq P 
\\
&[\bar{A}_{\mathrm{cl}}]_i^\top p_i [\bar{A}_{\mathrm{cl}}]_i\preceq \bar{P}_i,~i\in\mathbb{N}_1^{n_x} \\
&[\bar{C}_{\mathrm{cl},j}]_i^\top q_i [\bar{C}_{\mathrm{cl},j}]_i\preceq \bar{Q}_{i,j},~i\in\mathbb{N}_1^{n_x},~j\in\mathbb{N}_1^{p-1}.
\end{align*}
To obtain an LMI, we define new decision variables $X=P^{-1}$, $Y=KX$ that uniquely define $P,K$ with $X=\mathrm{diag}(x_1,\dots,x_{n_x})$. 
Furthermore, we define $\bar{P}_{x,i}=X\bar{P}_iX$, $\bar{Q}_{x,i,j}=X\bar{Q}_{i,j}X$. 
Multiplying the conditions by $X$ from left and right yields
\begin{align*}
\sum_{i=1}^{n_x} \bar{P}_{x,i}+ XQX+\sum_{j=1}^{p-1}\sum_{i=1}^{n_x} \bar{Q}_{x,i,j}+  Y^{\top}R_pY \preceq& X 
\\
([\bar{A}]_iX+[\bar{B}]_iY)^\top p_i ([\bar{A}]_iX+[\bar{B}]_iY)\preceq& \bar{P}_{x,i} \\
([\bar{C}_j]_iX+[\bar{D}_j]_iY)^\top q_i ([\bar{C}_j]_iX+[\bar{D}_j]_iY)\preceq &\bar{Q}_{x,i,j}.
\end{align*}
Finally, we use the Schur complement to ensure the inequalities are linear in the decision variables.
For the last two conditions, we obtain the equivalent LMIs:
\begin{align}
&\begin{bmatrix}
x_i&[\bar{A}]_iX+[\bar{B}]_iY\\
([\bar{A}]_iX+[\bar{B}]_iY)^\top & \bar{P}_{x,i}    
\end{bmatrix}\succeq 0 \label{eq:LMI_1}\\
&\begin{bmatrix}
1/q_i  & [\bar{C}_j]_iX+[\bar{D}_j]_iY\\
([\bar{C}_j]_iX+[\bar{D}_j]_iY)^\top&\bar{Q}_{x,i,j}
\end{bmatrix}\succeq 0.\label{eq:LMI_2}
\end{align}
We introduce the linear constraint
\begin{align}
&\overline{PQ}=\sum_{i=1}^{n_x}(\bar{P}_{x,i}+\sum_{j=1}^{p-1}\bar{Q}_{x,i,j}),\label{eq:LMI_4}
\end{align}
which allows us to write the first condition as 
\begin{align*}
\overline{PQ}\overline{PQ}^{-1}\overline{PQ}+ XQX+  Y^{\top}R_pY \preceq X.
\end{align*}
Applying the Schur complement yields
\begin{align}
&\begin{bmatrix}
X&\overline{PQ} &(R^{1/2}Y)^\top&(Q^{1/2}X)^\top \\
\overline{PQ}&\overline{PQ} &0&0\\
R^{1/2}Y&0&I&0\\
Q^{1/2}X&0&0&I
\end{bmatrix}\succeq 0 \label{eq:LMI_3}
\end{align}
%
Hence, Inequalities~\eqref{eq:LMI_1}, \eqref{eq:LMI_2} $\forall i\in\mathbb{N}_1^{n_x}$, $j\in\mathbb{N}_1^{p-1}$ and \eqref{eq:LMI_4}, \eqref{eq:LMI_3} ensure satisfaction of~\eqref{eq:terminal_cost}. 
These conditions are LMIs $X,Y,\bar{Q}_{x,i,j},\bar{P}_{x,i}$, $\overline{PQ}$, 
Notably, ensuring~\eqref{eq:LMI_1}--\eqref{eq:LMI_2} for all $\theta\in\Theta$ requires only the vertices of the $j$-step predictor for the $i$-th state, i.e., $2^{n_x+jn_u}$ vertices.  
Condition~\eqref{eq:LMI_4}--\eqref{eq:LMI_3} does not depend on $\theta$ and requires no vertices. 

\subsection{Additional details for the numerical example} \label{A:example}
A dataset of $(u,x)$ pairs was collected as in~\cite[Sec.~IV]{terzi2022robust} by exciting the system. Specifically, a random input signal $u$ was applied, with values in $\mathbb{R}_{-1}^{1}$. The system was also subject to a disturbance characterized by $|w_k|_\infty\leq 0.01$ applied to the state space-system with a matrix $M=[1,\ 0.1, \ 0.1]^\top$. This dataset was subsequently used to calculate the parameter set $\Theta$ as defined in Equ.~\eqref{eq:modelUncer}. It was also employed to derive estimates for the bounds $\mathbb{W}_x$ and $\mathbb{W}_y$, following the methodology outlined in \cite{lauricella2020set,terzi2022robust}.
The same dataset was also utilized to compute the worst-case prediction bound $\tau_{p}$, as defined in \cite[Equ.~(8)]{terzi2022robust}. This bound is necessary for the implementation of \cite{terzi2022robust}, where it is calculated for each predictor $j$ and each state $i$ as follows:
\begin{align}
\tau_{j,i} = &\max_{\varphi\in\mathbb{X}\times\mathbb{U}}{\max_{\theta \in \Theta_{j,i}}{\left| \varphi_{j}^\top (\theta-\hat{\theta}_{j,i}) \right|}}+\hat{\bar{\varepsilon}}{j,i},\\
&\forall j\in\mathbb{N}_1^p, \ i \in \mathbb{N}_1^3, \nonumber
\end{align}
where, $\hat{\bar{\varepsilon}}_{j,i}$ is estimated as described in \cite[Equ.~(14,15)]{lauricella2020set}. Based on this bound a minimal RPI set as been computed as shown in~\cite{terzi2022robust}. The nominal parameters vector $\hat{\theta}$ was selected to minimize the worst-case prediction error bound $\tau_p$, as outlined in \cite[Sect.~IV]{lauricella2020set}.
The auxiliary control law $K$ and the matrix $P$ were computed by solving the LMI detailed in Appendix~\ref{A:LMI}, while for~\cite{lorenzen2019robust} a similar procedure has been used considering model~\eqref{eq:system}.
Finally, we selected the low-complexity polytope $\pazocal{X}_0$ for both the proposed approach and the homothetic tube MPC~\cite{lorenzen2019robust} following the procedure described in~\cite[Sec.~5.4]{kouvaritakis2016model}.
\end{document}